\documentclass{amsart}
\usepackage{latexsym}
\usepackage{amsmath}
\usepackage{amssymb}
\begin{document}
\pagestyle{plain}
\title{Construction of the Thermodynamic Jamming Limit for the Parking Process and Other Exclusion Schemes on $\mathbb{Z}^{d}$}
\author{Thomas Logan Ritchie}
\address{Instituto de Matem\'{a}tica e Estat\'{\i}stica\\Universidade de S\~{a}o Paulo\\Rua do Mat\~{a}o 1010, CEP 05508-900, S\~{a}o Paulo SP, Brazil }
\email{thomasl@ime.usp.br}
\thanks{Sponsored by FAPESP: processo 00/09052-3}
\date{July 04, 2004}
\begin{abstract}
We provide an explicit construction of the thermodynamic jamming limit  for the parking process and other finite range exclusion schemes on $\mathbb{Z}^{d}$. By means thereof, a strong law of large numbers for occupation densities is accomplished, and, amongst other results, the so called ``super-exponential'' (i.e. gamma) decay of  pair correlation functions is established.
\end{abstract}
\maketitle

\section{Introduction}
\label{introduction}
In its simplest form the \emph{car-parking process} on $\mathbb{Z}^{d}$ (also called  \emph{random sequential adsorption} (RSA) in physical sciences) may be (informally) described as follows:\vspace{1mm}\\
\textit{ a site is chosen uniformly over the box $\Lambda_{n}^{(d)}=\{-n,\ldots,0,\ldots,n\}^{d}$ (the target region) and is occupied (spin/state 1) provided all its nearest neighbours are vacant (spin/state 0). Once occupied a site remains so for ever. The process continues until all sites in $\Lambda_{n}^{(d)}$ are either occupied or have at least one of their nearest neighbours occupied (the so called \emph{jamming limit}). Call $P_{n}$ the terminal configuration of $\Lambda_{n}^{(d)}$.}\vspace{1mm}\hspace{3mm}\footnote{Throughout this paper we shall abuse terminology and use $P_{n}$ to denote not only the terminal configuration/jamming limit, but also the process/algorithm leading to it. Since only  terminal configurations will be taken into account,  this identification is reasonable.}\\
This model as many of its generalizations (including longer range interactions and particles of different kinds) have long received attention in physical sciences, wherein issues like kinetics  near the jamming limit, spatial correlations, assymptotics (as the target region increases to infinity) of the first moments of $|P_{n}|$ (the ultimate number of adsopted particles, i.e. occupied sites), amongst others, have been studied by a series of non-rigorous theoretical methods, computer simulations and experiments~\cite{Evans,Talbot}.

Until recently, mathematical results have been restricted to one-\!dimensional models~\cite{Flory,Mackenzie,Keller, Renyi,Page,correlations} and to quasi-\!one-\!dimensional models~\cite{Fan,Baram1,Baram2}.

In 2002 M. D. Penrose~\cite{Penrose2} proved, amongst other results, a law of large numbers ($L^{p}$ -convergence) for the sequence $\left(|P_{n}|/|\Lambda_{n}^{(d)}|\right)_{n\geq0}$ for a large class of exclusion models (including those with finitely many spin states), treated therein as (coupled) finite interacting particle systems.

In this paper we construct in a common probability space (i.e. couple) a limit ``parking'' algorithm $P$, having the whole $\mathbb{Z}^{d}$ lattice as target region, together with the finite parking algorithms $P_{n}$, \small $n\geq0$ \normalsize in such a way that $P_{n}\rightarrow P$ almost surely. This construction being sufficiently general to cope with all sort of finite range (exclusion) interactions in a $\{0,1\}$ \emph{spin} state scenary. In this set up $P$ may be interpreted as the jamming limit of a virtual parking process over the whole $\mathbb{Z}^{d}$ lattice. We call $P$ the thermodynamic limit of the sequence $(P_{n})_{n\geq0}$, and its distribution $\mu$  the corresponding thermodynamic limit measure.

By means of this construction (section~\ref{construction}), we prove a strong law of large numbers for the sequence $\left(|P_{n}|/|\Lambda_{n}^{(d)}|\right)_{n\geq0}$, via the classical multiparameter ergodic theorem (section~\ref{SLLNs}); the super-exponential (i.e. gamma) decay of pair correlation functions (subsection~\ref{sedpcfs}); the super-exponential speed of convergence\footnote{In a sense to be specified in subsection~\ref{speed}.} to the thermodynamic limit measure $\mu$ (subsection~\ref{speed}) and discuss how to obtain lower and upper estimates for $\rho$, the almost sure limit of $\left(|P_{n}|/|\Lambda_{n}^{(d)}|\right)_{n\geq0}$ (section~\ref{conclusion}).

\section{Basic Definitions and Notation}
\label{basic}
The following structures will be basic blocks in the forthcoming constructions. Therefore we take the opportunity to establish the notation to be used henceforth.
\subsection*{ Probability and Topological Spaces}
The random objects we shall deal with in this paper will be defined on a fundamental probability space (i) and assume values in a Polish measure  space\footnote{a measure space whose $\sigma$-algebra is generated by a topology which is metrizable, complete and separable.}  (ii). As construction is the realm of this paper, we make these spaces explicit: 
\begin{description}
\item[(i)]
$\left( (0,1)^{\mathbb{Z}^{d}},\mathcal{B}_{(0,1)}^{\mathbb{Z}^{d}},
\lambda_{(0,1)}^{\mathbb{Z}^{d}}      \right)$, where $\mathcal{B}_{(0,1)}$ is the usual Borel $\sigma$-algebra in $(0,1)$, $\lambda_{(0,1)}$ is the Lebesgue (uniform) measure in $(0,1)$  and $\cdot^{\mathbb{Z}^{d}}$ means product in the appropiate sense;
\item[(ii)]
$\left( \{0,1\}^{\mathbb{Z}^{d}},\mathcal{B}_{\{0,1\}}^{\mathbb{Z}^{d}}\right)$, where 
$\mathcal{B}_{\{0,1\}}$ is taken for the power set of $\{0,1\}$,\\ $\mathcal{P}\left(\{0,1\} \right)\stackrel{\mathrm{def}}{=}\left\{\{0,1\},\{0\},\{1\},\emptyset\right\}$.                   
\end{description} 

Since both $(0,1)$ and $\{0,1\}$ are Polish topological spaces (under the euclidean and discrete metric respectively), so are $(0,1)^{\mathbb{Z}^{d}}$ and 
$\{0,1\}^{\mathbb{Z}^{d}}$ (under the respective product topologies). Moreover the product $\sigma$-algebras $\mathcal{B}_{(0,1)}^{\mathbb{Z}^{d}}$ and $\mathcal{B}_{\{0,1\}}^{\mathbb{Z}^{d}}$ are the respective Borel $\sigma$-algebras of these Polish spaces. Therefore the measure spaces above, viz
$\left( (0,1)^{\mathbb{Z}^{d}},\mathcal{B}_{(0,1)}^{\mathbb{Z}^{d}} \right)$ and
$\left( \{0,1\}^{\mathbb{Z}^{d}},\mathcal{B}_{\{0,1\}}^{\mathbb{Z}^{d}}\right)$, are naturally related to the Polish spaces $(0,1)^{\mathbb{Z}^{d}}$ and $\{0,1\}^{\mathbb{Z}^{d}}$ - product topologies implicit- in the sense that the $\sigma$-algebras of the former are the Borel $\sigma$-algebras of the latter.

These remarks may seem pedantic at first glance, but they are the key to conclude convergence of distribution from (almost sure) pointwise convergence of general random objects, as we shall see below.

\subsection*{Graphs}

In what follows the notion of vicinity between sites of the $\mathbb{Z}^d$ lattice will be implemented by the graph $\mathcal{G}^{(d)}_{\nu}\stackrel{\mathrm{def}}{=}\left( \mathbb{Z}^{d},\mathcal{E}^{(d)}_{\nu}   \right)$, where $\mathcal{E}^{(d)}_{\nu} $ is taken for set of bonds joining distinct sites $x,y$ in $\mathbb{Z}^{d}$, such that $|x-y|_{\sup}\stackrel{\mathrm{def}}{=}\max_{1\leq i\leq d}\left| x_{i}-y_{i}\right|\leq \nu$, for some distance $\nu \in \mathbb{N}^{*}\stackrel{\mathrm{def}}{=}\mathbb{N}\setminus \{0\}$.

\subsection*{Parking Schemes}
A subset $S^{(d)}_{\nu}\subset\{ 0,1 \}^{\Lambda^{(d)}_{\nu}}$, where   $\Lambda^{(d)}_{\nu}\stackrel{\mathrm{def}}{=}\{-\nu,\ldots,0,\ldots,\nu\}^{d}$ is the box/ball with centre at the origin and radius $\nu$ in the $\mathbb{Z}^{d}$ lattice, such that 
$ S^{(d)}_{\nu}\ni 0^{\Lambda^{(d)}_{\nu}}$ will be called a \emph{parking scheme} (with radius $\nu$ of interaction) on the $\mathbb{Z}^{d}$ lattice.\footnote{As it will become clear from the forthcoming definitions, we may suppose, without loss of generality, that $\eta(0^{d})=0$ for all $\eta\in S_{\nu}^{(d)}$.\\  Henceforth we shall abuse notation and write $0$ for $0^{d}$, whenever no confusion may arise.}

\subsection*{Parking Processes}
For a given \emph{parking scheme} $S^{(d)}_{\nu}$, we define the \emph{parking process} on the box $\Lambda^{(d)}_{n}$ with null boundary condition, $\left( \Lambda^{(d)}_{n}\ S^{(d)}_{\nu}\right)PP\mathit{nbc}$, to be the following (random) algorithm:
\begin{description}
\item[Step 1]set $P_{n}=0^{\mathbb{Z}^{d}}$;
\item[Step 2]choose a site $x\in\Lambda^{(d)}_{n}$ uniformly amongst the sites not chosen previously;
\item[Step 3]if $\left. \theta_{x}\left(P_{n}\right)\right|\Lambda^{(d)}_{\nu}\in S^{(d)}_{\nu}$, then $P_{n}(x)\leftarrow 1$;
\item[Step 4]if there are  points in $\Lambda^{(d)}_{n}$ not chosen yet, then go back to step 2, else stop the algorithm.
\end{description}

Above, for all $y\in\mathbb{Z}^{d}$, $\theta_{y}$ is the translation operator, i.e. 
for all $P\in\{0,1\}^{\mathbb{Z}^{d}}$, $\theta_{y}\left(P\right)(x)\stackrel{\mathrm{def}}{=}P(x+y)$. And                
$\left.P\right|\Lambda^{(d)}_{\nu}$ means the restriction of $P$ to the box $\Lambda^{(d)}_{\nu}$, that is, $\left.P\right|\Lambda^{(d)}_{\nu}\in \{0,1\}^{\Lambda^{(d)}_{\nu}}$ and for all $x\in\Lambda^{(d)}_{\nu}$, $\left.P\right|\Lambda^{(d)}_{\nu}(x)\stackrel{\mathrm{def}}{=}P(x)$.

It should be clear that the algorithm described above generates a random element $P_{n}$ in the Polish measure space $\left( \{0,1\}^{\mathbb{Z}^{d}},\mathcal{B}_{\{0,1\}}^{\mathbb{Z}^{d}}\right)$ and therefore a corresponding probability measure is induced therein.

We put these ideas into a \emph{probability setting} by means of
\newtheorem{definition}{Definition}[section]
\begin{definition}
\label{parkingprocessonboxes}
For a given \emph{parking scheme} $S^{(d)}_{\nu}$, we define the parking process on the box $\Lambda^{(d)}_{n}$ with null boundary condition, $\left( \Lambda^{(d)}_{n}\ S^{(d)}_{\nu}\right)PP\mathit{nbc}$, to be the following  algorithm:

For each $\omega \in (0,1)^{\mathbb{Z}^{d}}$,
\begin{description}
\item[Step 1]set $P_{n}(\omega)=0^{\mathbb{Z}^{d}}$;
\item[Step 2]choose  $x\in\Lambda^{(d)}_{n}$  such that \\ $\omega(x)=\inf\{\omega(\xi):\xi \in \Lambda^{(d)}_{n} \mbox{and $\xi$ has not been chosen previously}\}$;
\item[Step 3]if $\left. \theta_{x}\left(P_{n}(\omega)\right)\right|\Lambda^{(d)}_{\nu}\in S^{(d)}_{\nu}$, then $P_{n}(\omega)(x)\leftarrow 1$;
\item[Step 4]if there are  points in $\Lambda^{(d)}_{n}$ not chosen yet, then go back to step 2, else stop the algorithm.
\end{description}

\end{definition}
\hfill $\Box$\\
The formal algorithm described above is probabilisticaly equivalent to the former (informal) algorithm in the sense that both induce the same probability measure $\mu_{n}$ in the Polish measure space $\left( \{0,1\}^{
\mathbb{Z}^{d}},\mathcal{B}_{\{0,1\}}^{\mathbb{Z}^{d}}\right)$. This equivalence being due to the fact that the minimum of an \emph{iid} finite family of uniform random variables is uniformly distributed amongst them.

We call the output $P_{n}(\omega)$ the \emph{jamming limit} of the parking process just described\footnote{In the forthcoming text we shall frequently abuse terminilogy and use $P_{n}$ to denote also the underlying parking process leading to it.} and observe that $P_{n}$ is a $\{0,1\}^{\mathbb{Z}^{d}}
$-valued random element on 
$\left( (0,1)^{\mathbb{Z}^{d}},\mathcal{B}_{(0,1)}^{\mathbb{Z}^{d}}\right.$
$,\left.\lambda_{(0,1)}^{\mathbb{Z}^{d}}      \right)$.

Finally, we remark that each $P_{n},n\geq0$, is well (i.e. uniquely) defined, except for a set of null probability; ambiguity might occur at $\omega$'s such that $\omega(x)=\omega(y)$ for at least two distinct sites $x,y \in \mathbb{Z}^{d}$.
We shall return to this point later on, but, as it usually happens in probability theory, that should pose us no difficulty at all.

A natural question that arises at this point is the existence (or not) of a thermodynamic limit probability measure $\mu$ for the sequence $(\mu_{n})_{n\geq1}$ of probability measures on  $\left( \{0,1\}^{\mathbb{Z}^{d}},\mathcal{B}_{\{0,1\}}^{\mathbb{Z}^{d}}\right)$ induced by the sequence $\left(P_{n}\right)_{n\geq1}$ of random elements just defined.
We shall delve into this question in the next section. 
\section{Construction of the Thermodynamic Limit}
\label{construction}
In what follows we show that, for a given parking scheme $S_{\nu}^{(d)}$, there exists a unique probability measure $\mu$ on $\left( \{0,1\}^{\mathbb{Z}^{d}},\mathcal{B}_{\{0,1\}}^{\mathbb{Z}^{d}}\right)$ such that $\mu_{n}\Rightarrow \mu$. \footnote{the symbol ``$\Rightarrow$'' denotes convergence in distribution.} 
We shall accomplish this task by means of constructing a random element/algorithm $P:\left( (0,1)^{\mathbb{Z}^{d}},\mathcal{B}_{(0,1)}^{\mathbb{Z}^{d}},
\lambda_{(0,1)}^{\mathbb{Z}^{d}}      \right) \longrightarrow \left( \{0,1\}^{\mathbb{Z}^{d}},\mathcal{B}_{\{0,1\}}^{\mathbb{Z}^{d}}\right)$ such that $P_{n}\rightarrow P \ a.s.$\ \footnote{pointwise convergence regarding the product topology in $\{0,1\}^{\mathbb{Z}^{d}}$.}

For that purpose we set the following definitions:
\begin{definition}
A path of points in the graph $\mathcal{G}_{\nu}^{(d)}$ is a finite sequence $\left(x_{i}\right)_{0\leq i \leq n}$, $x_{i}\in \mathbb{Z}^{d}$; such that $\left|x_{i+1}-x_{i}\right|_{\sup}\leq \nu$, $0\leq i \leq n-1$.\hfill  $\Box$
\end{definition}
In what follows we shall say that the path $(x_{i})_{0\leq i\leq n}$ has length $n$ (rather than $n+1$), bearing in mind that $n$ bonds participate in such a path.
\begin{definition}
A path $\left(x_{i}\right)_{0\leq i \leq n}$ is said to be decreasing (subject to $\omega$\small $\in$\normalsize $(0,1)^{\mathbb{Z}^{d}}$), whenever the numerical sequence $\left(x_{i}(\omega)\right)_{0\leq i \leq n}$ is (strictly) decreasing.\footnote{$x_{i}(\omega)$ is taken for $\omega(x_{i})$, which should be formally prefered.} In this case we write $\downarrow_{(\omega)}\left(x_{i}\right)_{0\leq i \leq n}$, or simply $\downarrow\left(x_{i}\right)_{0\leq i \leq n}$, should $\omega$ be clear from the context.\hfill  $\Box$
\label{decreasingpath}
\end{definition}
It should be clear from definition~\ref{decreasingpath} that a decreasing path is always self-avoiding.
\begin{definition}
For $x,y\in\mathbb{Z}^{d}$ and $\omega\in(0,1)^{\mathbb{Z}^{d}}$, we say that $x$ is influenced by $y$ (or, equivalentely, that $y$ influences $x$) subject to $\omega$, if there exists a decreasing path $\downarrow_{(\omega)}\left(x_{i}\right)_{0\leq i \leq n}$, such that $x_{0}=x$ and $x_{n}=y$. In this case we write $y\uparrow_{(\omega)}x$, or simply $y\uparrow x$, should $\omega$ be clear from the context.\hfill  $\Box$
\end{definition}

We remark that the $\uparrow$ relation is both transitive and reflexive. 
\begin{definition}
\label{armour}
Given a  subset $X\subset\mathbb{Z}^{d}$, we define the \emph{armour} of $X$ (subject to $\omega\in(0,1)^{\mathbb{Z}^{d}}$) to be the (random) subset
\begin{equation*}
\mathcal{A}(X)(\omega)\stackrel{\mathrm{def}}{=}
\left\{y\in\mathbb{Z}^{d}:\exists x\in X:y\uparrow_{(\omega)}x \right\}  
\end{equation*}
\hfill  $\Box$
\end{definition}
As before, we shall omit $\omega$ whenever no confusion may arise.\\
It is clear from definition~\ref{armour} that, for $X,Y\subset\mathbb{Z}^{d}$ and $\omega\in(0,1)^{\mathbb{Z}^{d}}$,
\begin{description}
\item[(i)]$\mathcal{A}(X)(\omega)\supset X$;
\item[(ii)]$X=\biguplus_{i=1}^{n}X_{i}\Rightarrow\mathcal{A}(X)(\omega)=\bigcup_{i=1}^{n}\mathcal{A}(X_{i})(\omega)$;  \ \footnote{$\biguplus$ denotes disjoint union.}
\item[(iii)]$X\subset Y\Rightarrow\mathcal{A}(X)(\omega)\subset\mathcal{A}(Y)(\omega)$.
\end{description}
What may not be so clear is the following lemma.
\newtheorem{lemma}{Lemma}[section]
\begin{lemma}
For any finite $X\subset\mathbb{Z}^{d}$, $\mathcal{A}(X)$ is (almost surely) finite.
\label{finite}
\end{lemma}
\textbf{Proof:}\\
We divide the proof into three parts:\vspace{3mm}\\
\textbf{Part (i)}\vspace{1.5mm}\\
Suppose initially that $X=\{0\}$ and consider the event $A_{n}\stackrel{\mathrm{def}}{=}\left(\mathcal{A}(X)\not\subset\Lambda_{n\nu}^{(d)} \right)$, $n\in\mathbb{N}$.
 $A_{n}$ amounts to saying that there exists $y\in\mathbb{Z}^{d}\setminus\Lambda_{n\nu}^{(d)}$ such that $y\uparrow 0$, that is, there exists a decreasing (and hence self-avoiding) path $(x_{i})_{0\leq i\leq m}$, $m> n$, starting at 0 and ending outside the box $\Lambda_{n\nu}^{(d)}$. Thus $A_{n}$ occurs only if there exists a decreasing path of length $n$, $(x_{i})_{0\leq i\leq n}$ starting at $0$.
  Now observe that the probability that an arbitrary (but fixed) self-avoiding path $(x_{i})_{0\leq i\leq n}$ is decreasing is $1/(n+1)!$ and that the total number of self-avoinding  paths starting at $0$ with length $n$ is no larger than $(2\nu+1)^{d.n}$, 
to conclude (by subaditivity) that 
\begin{equation}
\label{estimation}
\mathbb{P}(A_{n})\leq \frac{(2\nu+1)^{d\cdot n}}{(n+1)!}\leq  \frac{(2\nu+1)^{d\cdot n}}{n!}
\end{equation}
Since the sequence $(A_{n})_{n\geq 0}$ is decreasing and $\mathbb{P}(A_{n})\rightarrow 0$ by (\!\!~\ref{estimation}), we must have $\mathbb{P}\left(\bigcap_{n=0}^{\infty}A_{n}\right)=0$. Now observe that  $\left\{\#\mathcal{A}(X)=\infty\right\}=\bigcap_{n=0}^{\infty}A_{n}$ to conclude that $\mathbb{P}\left(\left\{\#\mathcal{A}(X)<\infty\right\}=\left\{\#\mathcal{A}(X)=\infty\right\}^{c}\right)=1$, that is $\mathcal{A}\left(\{0\}\right)$ is almost surely finite.\vspace{3mm}\\
\textbf{Part (ii)}\vspace{1.5mm}\\
Set $X=\{x\}$, for some $x\in\mathbb{Z}^{d}$.\\
In this case we observe that $\left\{\#\mathcal{A}(X)<\infty\right\}=\theta_{-x}\left(\left\{\#\mathcal{A}\left(\{0\}\right)<\infty\right\}\right)$ and make use of the stationarity of the product measure $\lambda_{(0,1)}^{\mathbb{Z}^{d}}$ to conclude that $\mathbb{P}\left(\left\{\#\mathcal{A}(X)<\infty\right\}\right)$ $=$ $\mathbb{P}\left(\left\{\#\mathcal{A}\left(\{0\}\right)<\infty\right\}\right)=1$, by part (i) above, that is $\mathcal{A}\left(\{x\}\right)$ is almost surely finite for all $x\in\mathbb{Z}^{d}$\vspace{2mm}.\\
\textbf{Remark:} we might have observed the event equality \[\mathcal{A}\left(\{x\}\right)\not\subset B(x,n\nu)=\theta_{-x}\left(\mathcal{A}\left(\{0\}\right)\not\subset\Lambda_{n\nu}^{(d)}\right)\] made use of  stationarity  to yield \[\mathbb{P}\left(\mathcal{A}\left(\{x\}\right)\not\subset B(x,n\nu)\right)=\mathbb{P}\left(\theta_{-x}\left(\mathcal{A}\left(\{0\}\right)\not\subset\Lambda_{n\nu}^{(d)}\right)\right)=\mathbb{P}\left(\mathcal{A}\left(\{0\}\right)\not\subset\Lambda_{n\nu}^{(d)}\right)\]
and followed the same steps of part (i) to see that $\mathcal{A}\left(\{x\}\right)$ is almost surely finite.

This \emph{translation technique} will be recurrent in the forthcoming proofs.\vspace{3mm} \\
\textbf{Part (iii)}\vspace{1.5mm}\\
Suppose $X$ is a finite set, that is $X=\biguplus_{i=1}^{n}\{x_{i}:x_{i}\in X\}$. Thus by the remark after definition~\ref{armour} $\mathcal{A}(X)=\bigcup_{i=1}^{n}\mathcal{A}\left(\{x_{i}\}\right) \Rightarrow\left\{\#\mathcal{A}(X)=\infty\right\}=\bigcup_{i=1}^{n}\left\{\#\mathcal{A}\left(\{x_{i}\}\right)=\infty\right\}$ $\Rightarrow \mathbb{P}\left(\left\{\#\mathcal{A}(X)=\infty\right\}\right)\leq\sum_{i=1}^{n}\mathbb{P}\left(\left\{\#\mathcal{A}\left(\{x_{i}\}\right)=\infty\right\}\right)=0$ by part (ii) above. That is $\mathcal{A}(X)$ is almost surely finite.
\begin{flushright}$\blacksquare$\end{flushright}
In what follows it will be convinient to generalize \emph{parking processes} on boxes (definition~\ref{parkingprocessonboxes}) to parking processes on arbitrary finite subsets of $\mathbb{Z}^{d}$. This generalization does not involve any difficulty and is performed in the natural way.

\begin{definition}
\label{ppogs}
Let $X$ be a finite subset of $\mathbb{Z}^{d}$. For a given \emph{parking scheme} $S^{(d)}_{\nu}$, we define the parking process on $X$ with null boundary condition, \\
$\left(X\  S^{(d)}_{\nu}\right)PP\mathit{nbc}$, to be the following  algorithm:

For each $\omega \in (0,1)^{\mathbb{Z}^{d}}$,
\begin{description}
\item[Step 1]set $P_{X}(\omega)=0^{\mathbb{Z}^{d}}$;
\item[Step 2]choose  $x\in X$  such that \\ $\omega(x)=\inf\{\omega(\xi):\xi \in X \mbox{and $\xi$ has not been chosen previously}\}$;
\item[Step 3]if $\left. \theta_{x}\left(P_{X}(\omega)\right)\right|\Lambda^{(d)}_{\nu}\in S^{(d)}_{\nu}$, then $P_{X}(\omega)(x)\leftarrow 1$;
\item[Step 4]if there are  points in $X$ not chosen yet, then go back to step 2, else stop the algorithm.

\end{description}
\hfill  $\Box$
\end{definition}
As before, we call the output $P_{X}(\omega)$ the \emph{jamming limit} of the parking process just defined and observe that $P_{X}$ is a $\{0,1\}^{\mathbb{Z}^{d}}
$-valued random element on $\left( (0,1)^{\mathbb{Z}^{d}},\mathcal{B}_{(0,1)}^{\mathbb{Z}^{d}}\right.$
$,\left.\lambda_{(0,1)}^{\mathbb{Z}^{d}}      \right)$.

Now we are ready for the definition of the limit random element/algorithm $P$:
\begin{definition}
\label{thermodynamiclimit}
The limit random element/algorithm \\$P:\left( (0,1)^{\mathbb{Z}^{d}},\mathcal{B}_{(0,1)}^{\mathbb{Z}^{d}},
\lambda_{(0,1)}^{\mathbb{Z}^{d}}      \right)\rightarrow \left( \{0,1\}^{\mathbb{Z}^{d}},\mathcal{B}_{\{0,1\}}^{\mathbb{Z}^{d}}\right)$ is defined by  
\[P(\omega)(x)\stackrel{\mathrm{def}}{=}P_{\mathcal{A}(\{x\})(\omega)}(\omega)(x),\ \forall x\in \mathbb{Z}^{d}\]  

\ \hfill $\Box$
\end{definition}
In words definition~\ref{thermodynamiclimit} above means that for each $\omega\in(0,1)^{\mathbb{Z}^{d}}$ and $x\in \mathbb{Z}^{d}$, we determine the \emph{armour} of the set $\{x\}$, $\mathcal{A}\left(\{x\}\right)$, by means of definition~\ref{armour}, construct/implement the parking process over $\mathcal{A}\left(\{x\}\right)$ according to definition~\ref{ppogs} and then assign the value of the corresponding \emph{jamming limit} $P_{\mathcal{A}\left(\{x\}\right)}$ at site $x$ to $P(\omega)(x)$.
 Of course a parking scheme $S^{(d)}_{\nu}$ is implicit in definition~\ref{thermodynamiclimit}.
 
Our next step is to prove that $P_{n}\rightarrow P\  a.s.$ For that purpose we shall need the following lemma:
\begin{lemma}

$\left.P_{n}\right|\mathcal{A}\left(\{x\}\right)=
\left.P_{\mathcal{A} \left(\{x\}\right)}\right|\mathcal{A}\left(\{x\}\right)$ \emph{in} 
$\left\{\Lambda_{n}\supset\mathcal{A}\left(\{x\}\right)\right\}$.

\label{inclusion}
\end{lemma}
In words lemma~\ref{inclusion} asserts that the (random) sets $\mathcal{A}\left(\{x\}\right)$ are stable with respect to the sequence  $\left(P_{n}\right)_{n\geq 1}$; in the sense that, once $n$ becomes sufficiently large so that $\Lambda_{n}\supset\mathcal{A}\left(\{x\}\right)$, the \emph{jamming limits} $P_{m},m\geq n$ are identical over   
$\mathcal{A}\left(\{x\}\right)$.

 In the proof of lemma~\ref{inclusion} the notions of inner and outer boundaries of a set $X\subset \mathbb{Z}^{d}$ will be employed. We make them clear by means of the next definition.
 \begin{definition}
 For each $\nu\geq1$, the $\nu-$inner boundary of $X\subset \mathbb{Z}^{d}$ is defined to be the set
 \[\delta_{\nu}^{(in)}(X)\stackrel{\mathrm{def}}{=}\left\{x\in\mathbb{Z}^{d}:x\in X \mbox{ and }\exists y\in \mathbb{Z}^{d}\setminus X\mbox{ such that } |x-y|_{\sup}\leq\nu\right\}\]
 In the same fashion, the $\nu$-outer boundary of $X\subset \mathbb{Z}^{d}$ is defined to be the set
  \[\delta_{\nu}^{(out)}(X)\stackrel{\mathrm{def}}{=}\left\{x\in\mathbb{Z}^{d}:x\notin X \mbox{ and }\exists y\in X\mbox{ such that } |x-y|_{\sup}\leq\nu\right\}\]
 \hfill $\Box$
 \label{boundaries}
 \end{definition}
 For simplicity the subscrits will be suppressed, whenever $\nu$ is understood from the context.\vspace{1.5mm}
 
 \textbf{Remark:} an immediate consequence of definitions~\ref{armour} and~\ref{boundaries} is that $x(\omega)\leq y(\omega)$, whenever $x\in\delta^{(in)}\left(\mathcal{A}(X)\right)$, $y\in\delta^{(out)}\left(\mathcal{A}(X)\right)$ and  $|x-y|_{\sup}\leq\nu$, for all $X\subset\mathbb{Z}^{d}$. That is, an \emph{armour} has the noticeable property that in the interface with the outer enviroment its sites are numerically smaller than the \emph{interacting} outer sites.\vspace{3mm}
 
 \textbf{Proof of lemma~\ref{inclusion}}:
 Initially we observe that, for $x\in\mathbb{Z}^{d}$, (almost all) $\omega\in(0,1)^{\mathbb{Z}^{d}}$ induces a \emph{natural} order $\preccurlyeq$ in the set $\mathcal{A}\left(\{x\}\right)(\omega)\subset\mathbb{Z}^{d}$, in the sense that, for $y,z\in\mathcal{A}\left(\{x\}\right)$, $y\preccurlyeq z$, whenever $y(\omega)\leq z(\omega)$.
 
 Thus, for each $\omega$, one can construct a sequence $(y_{i})_{1\leq i\leq \#\mathcal{A}\left(\{x\}\right)(\omega)}$ in $\mathcal{A}\left(\{x\}\right)$ such that $y_{j}\preccurlyeq y_{k}$, whenever $j\leq k$; this sequence being (almost surely) uniquely determined.
 
 Now for a fixed $\omega$ and a \emph{time} $\tau\in\{1,2,\ldots,\#\mathcal{A}\left(\{x\}\right)(\omega)\}$, we observe the \emph{parking processes} $P_{n}$
 (definition~\ref{parkingprocessonboxes}) and $P_{\mathcal{A}\left(\{x\}\right)(\omega)}$ (definition~\ref{ppogs}) at time $\tau$, i.e. immediately before an \emph{arrival} (step 2 of the corresponding algorithms) at site $y_{\tau}$.
 
 We conclude the proof by induction (over $\tau$):\\
 \textbf{Base of Induction:} \\
 At time $\tau=1$, $P_{n}|\mathcal{A}\left(\{x\}\right)=P_{\mathcal{A}\left(\{x\}\right)}|\mathcal{A}\left(\{x\}\right)=0^{\mathcal{A}\left(\{x\}\right)}$, by step 1 of respective definitions.\\
 \textbf{Inductive Hypothesis:}\\
 Suppose that $P_{n}|\mathcal{A}\left(\{x\}\right)=P_{\mathcal{A}\left(\{x\}\right)}|\mathcal{A}\left(\{x\}\right)$ at time $\tau\in\{1,2,\ldots,\#\mathcal{A}\left(\{x\}\right)\}$.\\
\textbf{Inductive Step:}\\
If $y_{\tau}\notin\delta^{(in)}_{\nu}\left(\mathcal{A}(\{x\})\right)$, by the \emph{inductive hypothesis}, $P_{n}$ and $P_{\mathcal{A}\left(\{x\}\right)}$ are identical at time $\tau$ in the $\nu$-neighbourhood of $y_{\tau}$, i.e. $B(y_{\tau},\nu)\stackrel{\mathrm{def}}{=}\{y\in\mathbb{Z}^{d}:|y-y_{\tau}|_{\sup}\leq\nu\}$, hence $\theta_{y_{\tau}}(P_{n})|\Lambda_{\nu}^{(d)}=\theta_{y_{\tau}}(P_{\mathcal{A}(\{x\})})|\Lambda_{\nu}^{(d)}$.  Thus, by step 3 of the corresponding definitions, we see that $P_{n}|\mathcal{A}\left(\{x\}\right)=P_{\mathcal{A}\left(\{x\}\right)}|\mathcal{A}\left(\{x\}\right)$ at time $\tau+1$.\\
 If $y_{\tau}\in\delta^{(in)}_{\nu}\left(\mathcal{A}(\{x\})\right)$, again by the \emph{inductive hypothesis}, $P_{n}$ and $P_{\mathcal{A}\left(\{x\}\right)}$ coincide in $B(y_{\tau},\nu)\cap \mathcal{A}(\{x\})$. But, since $B(y_{\tau},\nu)\cap(\mathcal{A}(\{x\}))^{c}\subset\delta_{\nu}^{(out)}(\mathcal{A}(\{x\}))$, no $P_{n}$-arrival may have occured in this region up to time $\tau$ by the remark after definition~\ref{boundaries}, and therefore both $P_{n}$ and $P_{\mathcal{A}(\{x\})}$ must  equal $0$ thereon. Hence $P_{n}$ and $P_{\mathcal{A}(\{x\})}$ coincide again on the whole $B(y_{\tau},\nu)$ at time $\tau$, and by the same reason as before we must have  $P_{n}|\mathcal{A}\left(\{x\}\right)=P_{\mathcal{A}\left(\{x\}\right)}|\mathcal{A}\left(\{x\}\right)$ at time $\tau+1$.\hfill $\blacksquare$ 

Lemma~\ref{inclusion} is the key to the main result of this section:
\newtheorem{theorem}[lemma]{Theorem}
\begin{theorem}[Almost sure convergence to the limit algorithm $P$]
\label{convergence}
\[\mbox{\large $P_{n}\rightarrow P$ \normalsize\ \ \ } \lambda_{(0,1)}^{\mathbb{Z}^{d}}a.s.\]
\end{theorem}
\textbf{Remark:} before proving theorem~\ref{convergence} we emphasize that the \emph{a.s. convergence} therein is relative to the usual product topology  of  $\{0,1\}^{\mathbb{Z}^{d}}$.\vspace{3mm}\\
\textbf{Proof:}\\
Bearing in mind the foregoing remark, it is enough to show that for $\lambda_{(0,1)}^{\mathbb{Z}^{d}}$-almost all $\omega\in (0,1)^{\mathbb{Z}^{d}}$ and every $m\geq 0$, $P_{n}(\omega)|\Lambda_{m}^{(d)}=P(\omega)|\Lambda_{m}^{(d)}$, for $n$ sufficiently large.

Now, given $m\geq0$, take $n$ large enough so that $\mathcal{A}(\Lambda_{m}^{(d)})\subset\Lambda_{n}^{(d)}$. This is possible for almost all $\omega$ in view of lemma~\ref{finite}. For each $x\in\Lambda_{m}^{(d)}$, $\mathcal{A}\left(\{x\}\right)\subset\mathcal{A}(\Lambda_{m}^{(d)})$ (remark after definition~\ref{armour}). Hence 
\[  P(x)\stackrel{def.~\ref{thermodynamiclimit}}{=}P_{\mathcal{A}\left(\{x\}\right)}(x)\stackrel{lem.~\ref{inclusion}}{=}P_{n}(x),\ \forall x\in\Lambda_{m}^{(d)}\]
and we are done.\hfill \vspace{3mm}$\blacksquare$\\
Since the random elements $(P_{n})_{n\geq0}$ and $P$ take values in the Polish measure space $\left(\{0,1\}^{\mathbb{Z}^{d}},\mathcal{B}_{\{0,1\}}^{\mathbb{Z}^{d}}\right)$, almost sure convergence implies weak convergence of the corresponding distributions therein~\cite[prop.9 pag.355]{Gray}. Thus the following corollary:
\newtheorem{corollary}[lemma]{Corollary}
\begin{corollary}
If $(\mu_{n})_{n\geq0}$ and $\mu$ are the corresponding probability distributions of $(P_{n})_{n\geq0}$ and $P$ in $\left(\{0,1\}^{\mathbb{Z}^{d}},\mathcal{B}_{\{0,1\}}^{\mathbb{Z}^{d}}\right)$, then $\mu_{n}\Rightarrow\mu$. 
\end{corollary}

Before finishing this section we wish to emphasize that the construction of the limit measure carried out in definition~\ref{thermodynamiclimit} endows it with two noticeable properties, viz stationarity (i.e. translation invariace) and ergodicity.

In order to formalize this point a little further we observe that the random elements $(P_{n})_{n\geq0}$ and $P$ just defined (definitions~\ref{parkingprocessonboxes} and~\ref{thermodynamiclimit} respectively) are well defined, except for a set of probability zero, viz the set of $\omega$'s in $(0,1)^{\mathbb{Z}^{d}}$ that assign a same numerical value to at least two distinct sites of $\mathbb{Z}^{d}$ or on which an infinite descending path occurs.

This difficulty may be easily removed by setting $P_{n}=P=0^{\mathbb{Z}^{d}},\ n\geq0$, for instance, over this troublesome set. Once this is done, we can define the real measurable function $f$ on $\left( (0,1)^{\mathbb{Z}^{d}},\mathcal{B}_{(0,1)}^{\mathbb{Z}^{d}},\lambda_{(0,1)}^{\mathbb{Z}^{d}}      \right)$ by $f(\omega)\stackrel{\mathrel{def}}{=}P(\omega)(0)$ and observe that $P(\omega)(x)=f\left(\theta_{x}(\omega)\right)$, for all $x$ in $\mathbb{Z}^{d}$. Since the product measure $\lambda_{(0,1)}^{\mathbb{Z}^{d}}$ is ergodic and stationary in $\left( (0,1)^{\mathbb{Z}^{d}},\mathcal{B}_{(0,1)}^{\mathbb{Z}^{d}}\right)$, so will be $\mu$, the distribution of $P$, in $\left( \{0,1\}^{\mathbb{Z}^{d}},\mathcal{B}_{\{0,1\}}^{\mathbb{Z}^{d}}\right)$. This is a standard argument in ergodic theory and we refer the reader to ~\cite[chapt.6]{Durrett} for further details.

Therefore we have established another corollary to theorem~\ref{convergence}:
\begin{corollary}
\label{stationarity}
The thermodynamic limit measure $\mu$ is both stationary and ergodic in $\left( \{0,1\}^{\mathbb{Z}^{d}},\mathcal{B}_{\{0,1\}}^{\mathbb{Z}^{d}}\right)$.
\end{corollary}

The construction of the thermodynamic limit  measure $\mu$ carried out in this section corresponds to the concept of \emph{perfect simulation} in the sense of Fernandez, Ferrari and Garcia~\cite{Pablo}. In this technique, with the aid of a finite (yet random) number of i.i.d. (uniform) random variables, one generates (i.e. simulates) according to a prescribed thermodynamic measure a finite portion (window) of a infinite random object (field).

Besides the straightforward advantages comparing to ordinary \emph{Monte Carlo} simulations, wherein only finite approximations of the thermodynamic measure are actually simulated, this constructive technique is also very useful regarding theoretical proofs as the forthcoming sections testify.

\section{A Strong Law for Occupation Densities}
\label{SLLNs}
First of all we formalize the concept of occupation density refered to in section~\ref{introduction}:
\begin{definition}
For $n\in\mathbb{N}$, we define the \emph{occupation density} at the \emph{jamming limit} $P_{n}$ to be the random variable
\[\rho_{n}(\omega)\stackrel{\mathrel{def}}{=}\frac{\sum_{x\in\Lambda_{n}^{(d)}}P_{n}(\omega)(x)}{|\Lambda_{n}^{(d)}|=(2n+1)^{d}}\mbox{\hspace{7mm} on  $\left( (0,1)^{\mathbb{Z}^{d}},\mathcal{B}_{(0,1)}^{\mathbb{Z}^{d}},
\lambda_{(0,1)}^{\mathbb{Z}^{d}}      \right)$}\]
\hfill $\Box$
\label{sllns}
\end{definition}
As we mentioned in section~\ref{introduction}, the existence (or not) of $\lim_{n\rightarrow\infty}\mathbb{E}(\rho_{n})$ is a question of relevance in \emph{random sequential adsorption} (RSA) contexts. Untill recently an (affirmative) answer was known only for one-dimensional models~\cite{Flory,Keller,Mackenzie,Page,Renyi} or for quasi-one-dimensional models~\cite{Baram1,Baram2,Fan}. Recently Penrose~\cite{Penrose2} proved a law of large numbers (convergence in probability) for the sequence $(\rho_{n})_{n\geq1}$ for a large class of exclusion schemes (considered therein as finite particle systems) on $\mathbb{Z}^{d}$. This class overlaps to some extend to the class considered here.\footnote{We refer the reader to section~\ref{conclusion} for a discussion about this point.} Naturally, this result assigns an affirmative answer (by means of a straightforward application of the Bounded Convergence Theorem) to the existence of $\lim_{n\rightarrow\infty}\mathbb{E}(\rho_{n})$.

Here with the aid of the construction carried out in section~\ref{construction}, we shall accomplish a simple proof for a strong law of large numbers (SLLN) for  $(\rho_{n})_{n\geq1}$.

The first step to this proof is to observe that we can invoke \emph{Wiener's multiparameter ergodic theorem}
~\cite[Appendix 14.A]{Georgii} to conclude that

\begin{eqnarray}
\lim_{n\rightarrow\infty}\frac{\sum_{x\in\Lambda_{n}^{(d)}}P(\omega)(x)}{|\Lambda_{n}^{(d)}|}&=&
\lim_{n\rightarrow\infty}\frac{\sum_{x\in\Lambda_{n}^{(d)}}f\circ\theta_{x}(\omega)}{(2n+1)^{d}}
 =\mathbb{E}(f) \nonumber \\
\label{Wiener} &=&\mathbb{E}(f\circ\theta_{0})
=\mathbb{E}\left(P(0)\right)\stackrel{\mathrm{def}}{=}\rho\mbox{\hspace{10mm}}\lambda_{(0,1)}^{\mathbb{Z}^{d}}\mbox{-a.s.} 
\end{eqnarray}

\noindent
for some constant $\rho\in[0,1]$. This is true since $f$ is measurable and integrable on 
$\left( (0,1)^{\mathbb{Z}^{d}},\mathcal{B}_{(0,1)}^{\mathbb{Z}^{d}},
\lambda_{(0,1)}^{\mathbb{Z}^{d}}      \right)$ 
and the product measure is stationary and ergodic in $\left( (0,1)^{\mathbb{Z}^{d}},\mathcal{B}_{(0,1)}^{\mathbb{Z}^{d}}\right)$.

One might guess that the discrepancies between $P(\omega)|\Lambda_{n}^{(d)}$ and $P_{n}(\omega)|\Lambda_{n}^{(d)}$ are due to some kind of boundary effect. Next lemma quantifies the extension of this boundary pertubation.
\begin{lemma}
\label{pertubation}
If $A_{n}$ is the event $\left\{\mathcal{A}\left(\Lambda_{n}^{(d)}\right)\not\subset\Lambda_{(n+\lfloor\sqrt{n}\rfloor \nu)}^{(d)}\right\}$, then \\$\mathbb{P}(\limsup A_{n})=0$, i.e. with probability one the sequence $(A_{n})_{n\geq1}$ occurs finitely many times.\footnote{$\lfloor x\rfloor\stackrel{\mathrm{def}}{=}$ greatest integer less than or equal to $x$.}
\end{lemma}
\textbf{Proof:}\\
By means of the estimation (\!~\ref{estimation}) carried out in lemma~\ref{finite} and bearing in mind the stationarity of the product measure as we emphasized in the remark after part (ii) of lemma~\ref{finite},  subaditivity yields:
\begin{equation*}
\mathbb{P}(A_{n})\leq(2n+1)^{d}\cdot\frac{(2\nu+1)^{d\lfloor\sqrt{n}\rfloor}}{\lfloor\sqrt{n}\rfloor!}\ \ :=a_{n}
\end{equation*}
Therefore by Stirling's formula\footnote{$n!\simeq n^{n}e^{-n}\sqrt{2\pi n}$, i.e. $n^{n}e^{-n}\sqrt{2\pi n}/n!\rightarrow1$ as $n\rightarrow\infty$.}:
\begin{equation*}
\mathbb{P}(A_{n})\leq a_{n}\simeq\frac{1}{\sqrt{2\pi}}\cdot\frac{(2n+1)^{d}}{\lfloor\sqrt{n}\rfloor^{1/2}}\cdot\frac{(2\nu+1)^{d\lfloor\sqrt{n}\rfloor}}{\lfloor\sqrt{n}\rfloor^{\lfloor\sqrt{n}\rfloor}}\cdot e^{\lfloor\sqrt{n}\rfloor}\ \ :=b_{n}
\end{equation*}
But, since $b_{n}=o(1/n^{2})$,\ \footnote{Landau's notation: $f_{n}=o(g_{n})$ means that $f_{n}/g_{n}\rightarrow0$ as $n\rightarrow\infty$.} \[\sum_{n=1}^{\infty}b_{n}<\infty\Rightarrow\sum_{n=1}^{\infty}a_{n}<\infty\Rightarrow\sum_{n=1}^{\infty}\mathbb{P}(A_{n})<\infty.\]
\noindent
Thus the proof is concluded by an application of the first lemma of Borel-Cantelli.

\hfill \vspace{2mm}$\blacksquare$

Now we are ready for the main result of this section:
\begin{theorem}
\label{SLLN}
A Strong Law of Large Numbers for Occupation Densities:
\[\rho_{n}\rightarrow\rho\mbox{\hspace{10mm}}\lambda_{(0,1)}^{\mathbb{Z}^{d}}\mbox{-a.s.}\] 
\end{theorem}
\noindent\textbf{Proof:}\vspace{2mm}\\
 We know by lemma~\ref{pertubation} that the last occurrence of $(A_{n})_{n\geq1}$ is almost surely finite, that is, $\bar{n}=\sup \{n\in\mathbb{N}:A_{n} \makebox{ occurs}\}<\infty \ a.s.$ Given $\epsilon>0$, pick $i\in\mathbb{N}$ sufficiently large so that $i\in\left[n+\lfloor\sqrt{n}\rfloor\nu,(n+1)+\lfloor\sqrt{n+1}\rfloor\nu\right)$ for some $n>\bar{n}$ and such that $\left( \frac{2(m+1)+\lfloor\sqrt{m+1}\rfloor\cdot2\nu +1}{2m+1}\right)^{d}-1<\epsilon$, for all $m\geq n$. Therefore $P|\Lambda_{n}^{(d)}=P_{i}|\Lambda_{n}^{(d)}$ by the proof of theorem~\ref{convergence}. 
  
  Hence 

\begin{eqnarray*}
\left|\frac{\sum_{x\in\Lambda_{i}^{(d)}}P(x)-\sum_{x\in\Lambda_{i}^{(d)}}P_{i}(x)}{(2i+1)^{d}}\right|\leq &
\mbox{\hspace{0mm}}
\frac{\left\{\left[(n+1)+\lfloor\sqrt{n+1}\rfloor\nu\right]\cdot2+1\right\}^{d}-(2n+1)^{d}}{(2n+1)^{d}}
\\
=&\left(\frac{2(n+1)+\lfloor\sqrt{n+1}\rfloor\cdot2\nu+1}{(2n+1)} \right)^{d}-1<\epsilon
\end{eqnarray*}

That is 
\[\lim_{i\rightarrow\infty}\left(\frac{\sum_{x\in\Lambda_{i}^{(d)}}P(x)}{(2i+1)^{d}}-
\frac{\sum_{x\in\Lambda_{i}^{(d)}}P_{i}(x)}{(2i+1)^{d}}\right)=0 \mbox{\hspace{5mm}a.s.}\]

Since
 $\lim_{i\rightarrow\infty}\frac{\sum_{x\in\Lambda_{i}^{(d)}}P(x)}{(2i+1)^{d}}=\rho\ a.s.$ by(\!~\ref{Wiener}), $\rho_{i}=\frac{\sum_{x\in\Lambda_{i}^{(d)}}P_{i}(x)}{(2i+1)^{d}}\rightarrow\rho \ a.s.$ \hfill \vspace{2mm}$\blacksquare$

 \ \vspace{2mm}\\
  As we mentioned in the beginning of this section, a straightfoward application of the Bounded Convergence Theorem now garantees the existence of $\lim_{n\rightarrow\infty}\mathbb{E}(\rho_{n})$.
\begin{corollary}
$\lim_{n\rightarrow\infty}\mathbb{E}(\rho_{n})=\rho$\vspace{2mm}
\end{corollary}

In section~\ref{conclusion} we briefly discuss how to estimate the constant $\rho$ numerically.

\section{Super-exponential Behaviour}
\label{seb}
First of all we precise the meaning of \emph{super-exponential decay}.
\begin{definition}
\label{superex}
A function $f:\mathbb{N}\rightarrow\mathbb{R}$ is said to decay \emph{super-exponentially} whenever $f(n)=o(e^{-\alpha n}),\forall\alpha>0$, i.e. $\lim_{n\rightarrow \infty}f(n).e^{\alpha n}=0, \forall\alpha>0$.\hfill$\vspace{1 mm}\Box$
\end{definition}
In physicochemical literature the so called \emph{random sequential adsorptian} (RSA) processes are distinguished from \emph{equilibrium systems} by several noticeable properties~\cite{Evans,Talbot}. Amongst them we pay special attention to super-exponential decay of pair correlation functions, a feature of RSA models that differs remarkably from the (usual) exponential decay of equilibrium systems. Although described in physical and physicochemical literature~\cite{Evans,Talbot}, a rigorous account thereof is available only for one-dimensional models~\cite{correlations}.

Here super-exponential decay of pair-correlation functions is established from the construction of the thermodynamical limit measure $\mu$ carried out in section~\ref{construction}. This amounts to saying that this feature is a universal property of the whole class of exclusion schemes studied in this paper, regardless of dimension ($d$) and range of interaction ($\nu$).

Another universal property derived from section~\ref{construction} is the super-exponential \emph{speed} of convergence to the  limit measure $\mu$. We cope with this point in subsection~\ref{speed}.

\subsection{Super-Exponential Decay of Pair-correlation Functions}\label{sedpcfs}
   
We recall that the \emph{correlation coefficient} $\sigma_{X,Y}$ of two random variables $X,Y$ is defined by $\sigma_{X,Y}=\mathrm{Cov}(X,Y)/(\sigma_{X}\sigma_{Y})$, where $\mathrm{Cov}(X,Y)$ stands for the covariance of $X$ and $Y$, and $\sigma_{X}$, $\sigma_{Y}$ for the corresponding standard deviations.

Regarding the exclusion schemes studied here, one might guess how the state of site $x\in\mathbb{Z}^{d}$ influences the state of another site $y\in\mathbb{Z}^{d}$ with respect to the  limit measure $\mu$. Therefore we shall  investigate the assymptotic behaviour of the correlation coefficient $\sigma_{xy}$ of the random variables $P_{x}\stackrel{\mathrm{def}}{=}P(x)$ and $P_{y}\stackrel{\mathrm{def}}{=}P(y)$ as the distance $|x-y|_{sup}$ tends to infinity. Thanks to the translation invariace of $\mu$ (corollary~\ref{stationarity}) $ \sigma_{xy}=\sigma_{0(y-x)}$, and we may suppose $x=0$ without loss of generality.

The central result of this subsection (theorem~\ref{sedpcf}) asserts that $\sigma_{0x}$ decays super-exponentially as $x$ tends towards infinity. More precisely, this amounts to saying that $\lim_{x\rightarrow\infty}\sigma_{0x}e^{\alpha|x|}=0,\forall\alpha>0 $. In order to prove theorem~\ref{sedpcf}, it will be convenient to enlarge (i.e. square) our basic probability space $\left( (0,1)^{\mathbb{Z}^{d}},\mathcal{B}_{(0,1)}^{\mathbb{Z}^{d}},
\lambda_{(0,1)}^{\mathbb{Z}^{d}}      \right)$ to $\left( \left[(0,1)^{2}\right]^{\mathbb{Z}^{d}},\mathcal{B}_{(0,1)^{2}}^{\mathbb{Z}^{d}},
\lambda_{(0,1)^{2}}^{\mathbb{Z}^{d}}      \right)$. In this new setting each $\tilde{\omega}\in\left[(0,1)^{2}\right]^{\mathbb{Z}^{d}}$ may be thought of as an ordered pair $\tilde{\omega}=(\omega_{1},\omega_{2})$, for $\omega_{i}\in(0,1)^{\mathbb{Z}^{d}}$, $i=1,2$.

On this squared probability space we define the following projection operators taking values in the former unsquared measure space $\left( (0,1)^{\mathbb{Z}^{d}},\mathcal{B}_{(0,1)}^{\mathbb{Z}^{d}}\right)$:
\begin{definition}
\label{operators}
\ 
\begin{description}
\item[\textbf{Cardinal Projection Operators}]
\[\omega_{1}(\tilde{\omega})\stackrel{\mathrm{def}}{=}\omega_{1};\hspace{2mm}  \omega_{2}(\tilde{\omega})\stackrel{\mathrm{def}}{=}\omega_{2}.\hspace{74mm}\]
\item[Hybrid Projection Operators] 

\[\omega_{i,j}^{(n)}(\tilde{\omega})(x) \stackrel{\mathrm{def}}{=}\omega_{i}(\tilde{\omega})(x)1_{\{x_{1}\leq n\}}+\omega_{j}(\tilde{\omega})(x)1_{\{x_{1}> n\}};\ i,j\in\{1,2\},i\neq j, \]
for $n\in\mathbb{N}$ and bearing in mind the canonical representation $x=(x_{1},x_{2},$  $\ldots,x_{d})$, for sites in $\mathbb{Z}^{d}$ as d-tuples of integers.\footnote{$1_{\{x_{1}\leq n\}}=1$ if $x_{1}\leq n$, and $0$ otherwise \ldots}\hfill$\Box$
\end{description}
\end{definition}
In words definition~\ref{operators} denotes that $\omega_{i,j}^{(n)}$ plays the role of $\omega_{i}$  in the half hyper-space $\{x_{1}\leq n\}$ and the role of $\omega_{j}$ in $\{x_{1}> n\}$.

A straightfoward implication of definition~\ref{operators} is that all the projection operators above share the same probability distribution in $\mathcal{B}_{(0,1)}^{\mathbb{Z}^{d}}$, viz $\lambda_{(0,1)}^{\mathbb{Z}^{d}}$. Moreover, $\omega_{1}\perp\omega_{2}$ and $\omega_{1,2}^{(n)}\perp\omega_{2,1}^{(n)},\forall n\geq0$. \footnote{$\perp$ stands for independence.}

The foregoing paragraph amounts to saying that $P(\omega_{1})$, $P(\omega_{2})$, $P(\omega_{1,2}^{(n)})$, $P(\omega_{2,1}^{(n)})$ share the same probability distribution $\mu$ in $\mathcal{B}_{\{0,1\}}^{\mathbb{Z}^{d}}$ and that $P(\omega_{1,2}^{(n)})$ and $P(\omega_{2,1}^{(n)})$ are independent, i.e.
\begin{equation}
P(\omega_{1})\sim P(\omega_{2})\sim P(\omega_{1,2}^{(n)})\sim P(\omega_{2,1}^{(n)})\sim\mu,\forall n\geq0\mbox{\hspace{31mm}}
\label{identical}
\end{equation}
\begin{equation}
P(\omega_{1,2}^{(n)})\perp P(\omega_{2,1}^{(n)}),\forall n\geq0\mbox{\hspace{72mm}}
\label{independence}
\end{equation}
Now we are prepared to prove the key result of this subsection, viz\footnote{Henceforth, for  $x\in\mathbb{Z}^{d}$, $|x|$ stands for $|x|_{\sup}$.}
\begin{theorem}[Super-exponential decay of correlation functions]
\label{sedpcf}
\ \\
$\sigma_{0x}$ decays super-exponentially as $x\rightarrow\infty$, i.e. $\lim_{x\rightarrow\infty}\sigma_{0x}\cdot e^{\alpha|x|}=0,\forall\alpha>0$.
\end{theorem}
\textbf{Proof}:\\
Let $x=(x_{1},x_{2},\ldots,x_{d})\in\mathbb{Z}^{d}$ and suppose initially that $x_{1}=|x|=\sup_{1\leq i\leq d}\\\{|x_{i}|\}$. Set $y=\lfloor x_{1}/2\rfloor$ and pick the largest $n\in\mathbb{N}$ such that $2\nu n+4\nu\leq x_{1}$, i.e. $n=\lfloor x_{1}/2\nu\rfloor-2$. Hence\footnote{Remember that $\nu$ stands for the \emph{radius of interaction} of the underlying parking scheme $S_{\nu}^{(d)}$.} 
$n\nu+2\nu\leq y\leq|x|-(n\nu+2\nu)$.\\
Now observe that
$\mathcal{A}(\{0\})(\omega_{1})\subset B(0,n\nu)\Rightarrow\mathcal{A}(\{0\})(\omega_{1})=\mathcal{A}(\{0\})(\omega_{1,2}^{(y)})$. This is true because $\omega_{1}$ and $\omega_{1,2}^{(y)}$ are indistinguishable on the left of $y$  (definition~\ref{operators}). By the same token $\mathcal{A}(\{x\})(\omega_{1})\subset B(x,n\nu)\Rightarrow\mathcal{A}(\{x\})(\omega_{1})=\mathcal{A}(\{x\})(\omega_{2,1}^{(y)})$. Thus on the event $E\in\mathcal{B}_{(0,1)^{2}}^{\mathbb{Z}^{d}}$,
\[E\stackrel{\mathrm{def}}{=}\left(\mathcal{A}(\{0\})(\omega_{1})\subset B(0,n\nu)\right)\cap\left(\mathcal{A}(\{x\})(\omega_{1})\subset B(x,n\nu)\right)\]
we must have
\begin{equation*}
P(\omega_{1})(0)\stackrel{\mathrm{def.~\ref{thermodynamiclimit}}}{=}P_{\mathcal{A}(\{0\})(\omega_{1})}(\omega_{1})(0)=P_{\mathcal{A}(\{0\})(\omega_{1,2}^{(y)})}(\omega_{1,2}^{(y)})(0)\stackrel{\mathrm{def.~\ref{thermodynamiclimit}}}{=}P(\omega_{1,2}^{(y)})(0)
\end{equation*}

\begin{equation*}
P(\omega_{1})(x)\stackrel{\mathrm{def.~\ref{thermodynamiclimit}}}{=}P_{\mathcal{A}(\{x\})(\omega_{1})}(\omega_{1})(x)=P_{\mathcal{A}(\{x\})(\omega_{2,1}^{(y)})}(\omega_{2,1}^{(y)})(x)\stackrel{\mathrm{def.~\ref{thermodynamiclimit}}}{=}P(\omega_{2,1}^{(y)})(x)
\end{equation*}
and hence
\begin{equation}
\label{product}
P(\omega_{1})(0)\cdot P(\omega_{1})(x)=P(\omega_{1,2}^{(y)})(0)\cdot P(\omega_{2,1}^{(y)})(x)\mbox{\hspace{5mm} on $E$.\hspace{0mm}}
\end{equation}
Since 
\begin{equation}
\label{est}
\mathbb{P}(E^{c})\leq 2 \cdot \frac{(2\nu+1)^{d\cdot n}}{n!}
\end{equation}
by (\!~\ref{estimation}); (\!~\ref{product}) and (\!~\ref{est})  yield
\begin{equation}
\label{first}
\left|\mathbb{E}[P(\omega_{1})(0)\cdot P(\omega_{1})(x)]-\mathbb{E}[P(\omega_{1,2}^{(y)})(0) \cdot P(\omega_{2,1}^{(y)})(x)]\right| = 
\end{equation}

\[
\left|\mathbb{E}[P(\omega_{1})(0)\cdot P(\omega_{1})(x)-P(\omega_{1,2}^{(y)})(0)\cdot P(\omega_{2,1}^{(y)})(x)]\right| \leq \]

\[
\mathbb{E}\left(\left|([P(\omega_{1})(0)\cdot P(\omega_{1})(x)-P(\omega_{1,2}^{(y)})(0)\cdot P(\omega_{2,1}^{(y)})(x)\right|\right)\leq\]
\begin{equation*}  \mathbb{P}(E^{c})\leq 2\cdot \frac{(2\nu+1)^{d\cdot n}}{n!}\end{equation*}

\noindent
However
\begin{eqnarray}
\mathbb{E}[P(\omega_{1,2}^{(y)})(0)\cdot P(\omega_{2,1}^{(y)})(x)]&=&\mathbb{E}[P(\omega_{1,2}^{(y)})(0)]\cdot \mathbb{E}[P(\omega_{2,1}^{(y)})(x)] \nonumber \\ \label{second}
&=&\mathbb{E}[P(\omega_{1})(0)] \cdot \mathbb{E}[P(\omega_{1})(x)]
\end{eqnarray}
where we have used (\!\!~\ref{independence}) in the former equality and (\!\!~\ref{identical}) in the latter. Thus substituting (\!\!~\ref{second}) in (\!\!~\ref{first}) yield
\begin{equation}
\label{third}
\left|\mathbb{E}[P(\omega_{1})(0)\cdot P(\omega_{1})(x)]-\mathbb{E}[P(\omega_{1})(0)]\cdot \mathbb{E}[P(\omega_{1})(x)]\right|\leq  2\cdot \frac{(2\nu+1)^{d\cdot n}}{n!}
\end{equation}
Meanwhile, $P\circ\omega_{1}$ (on the enlarged/squared probability space) and $P$ (on the former/unsquared probability space) share the same distribution $\mu$ in $\mathcal{B}_{\{0,1\}}^{\mathbb{Z}^{d}}$ thanks to the remark after definition~\ref{operators}, i.e. $(P\circ\omega_{1})\sim P\sim\mu$. Therefore (\!\!~\ref{third}) amounts to
\begin{equation}
\left|\mathbb{E}(P_{0}\cdot P_{x})-\mathbb{E}(P_{0})\cdot \mathbb{E}(P_{x})\right|\leq  2\cdot \frac{(2\nu+1)^{d\cdot n}}{n!}
\end{equation}
Substituting $n$ in terms of $|x|$ as in the first paragraph of the current proof, and recalling the definition of covariance of two random variables, we finally have
\begin{equation}
\label{fourth}
\left|\mathrm{Cov}(P_{0},P_{x})\right|\leq 2\cdot \frac{(2\nu+1)^{d\cdot(\lfloor|x|/2\nu\rfloor-2)}}{(\lfloor|x|/2\nu\rfloor-2)!}
\end{equation}
At this point we observe that nothing would change in the foregoing calculations if it were the case that $x_{i}=\pm|x|$, $1\leq i\leq d$, in the beginning of this proof. A suitable redefinition of the projection operators $\omega_{1,2}^{(\cdot)}$ and $\omega_{2,1}^{(\cdot)}$ along the corresponding coordinate axes being all we would have to do to establish (\!\!~\ref{fourth}) once again. Therefore (\!\!~\ref{fourth}) remains true for all $x\in\mathbb{Z}^{d}$.

Having said this, since
\[\sigma_{0}^{2}\stackrel{\mathrm{def}}{=}\mathrm{Var}(P_{0})=\mathrm{Var}(P_{x})\stackrel{\mathrm{def}}{=}\sigma_{x}^{2},\ \forall x\in\mathbb{Z}^{d}\]
thanks to the stationarity of $\mu$ (corollary~\ref{stationarity}), we simply make use of the definition of correlation coefficient in (\!\!~\ref{fourth}) to establish that 
\begin{equation}
\label{fifth}
\left|\sigma_{0x}\right|\leq \frac{2}{\sigma_{0}^{2}}\cdot \frac{(2\nu+1)^{d\cdot(\lfloor|x|/2\nu\rfloor-2)}}{(\lfloor|x|/2\nu\rfloor-2)!}=o(e^{-\alpha|x|}),       \  \forall \alpha>0 
\end{equation}
and hence conclude the proof.\hfill$\blacksquare$

\subsection{Speed of Convergence to the Thermodynamic Limit}
\label{speed}
Since $P_{n}$ actually takes values in the finite set $\{0,1\}^{\Lambda_{n}^{(d)}}$, it is reasonable to quantify the rate at which the sequence $(\mu_{n})_{n\geq0}$ converges to the limit measure $\mu$ by estimating how fast the numerical sequence $|\mu(L)-\mu_{n}(L)|$ converges to zero, whenever $L$ is a \emph{local} set, i.e. a set determined by a finite number of sites in $\mathbb{Z}^{d}$.

The concept of \emph{locallity }can be formalized as follows:
\begin{definition}
\label{local}
A set $L\in \mathcal{B}_{\{0,1\}^{\mathbb{Z}^{d}}}$ is said to be local (or cylindrical) if for some $m\in\mathbb{N}$ there exists a subset $\tilde{L}\subset\{0,1\}^{\Lambda_{m}^{(d)}}$ such that \[\{0,1\}^{\mathbb{Z}^{d}}\ni\omega\in L\Leftrightarrow\{0,1\}^{\Lambda_{m}^{(d)}}\ni\omega|\Lambda_{m}^{(d)}\in\tilde{L}\] 
We write $\mathcal{L}$ for the class of all local sets in $\mathcal{B}_{\{0,1\}^{\mathbb{Z}^{d}}}$.\hfill$\Box$ 
\end{definition}
Next theorem asserts that the numerical convergence above is super-\hspace{.5mm}exponentially fast.
\begin{theorem}
If $L\in\mathcal{L}$, then $|\mu(L)-\mu_{n}(L)|=o(e^{-\alpha n}),\ \forall\alpha>0.$
\end{theorem}
\textbf{Proof}:\\ 
Choose $\tilde{L}$ and $m$ like in definition~\ref{local}. For $n>m$,
\begin{equation*}
\mathbb{P}\left( \mathcal{A}(\Lambda_{m}^{(d)})\not\subset\Lambda_{n}\right)\leq(2m+1)^{d}\cdot\frac{(2\nu+1)^{d\cdot\lceil\frac{n-m}{\nu}\rceil}}{\lceil\frac{n-m}{\nu}\rceil!}
\end{equation*}
by (\!\!~\ref{estimation}) and subaditivity.\footnote{$\lceil x\rceil\stackrel{\mathrm{def}}{=}$ the least integer greater than or equal to $x$.}

However $P|\Lambda_{m}^{(d)}=P_{n}|\Lambda_{m}^{(d)}$ on $\left(\mathcal{A}(\Lambda_{m}^{(d)})\subset\Lambda_{n}\right)$ by the proof of theorem~\ref{convergence}. Moreover
\begin{equation*}
\mu(L)=\mathbb{P}(\{P\in L\})=\mathbb{P}(\{P|\Lambda_{m}^{(d)}\in\tilde{L}\})=\mathbb{E}\left(1_{\{P|\Lambda_{m}^{(d)}\in\tilde{L}\}}\right)
\end{equation*}
and
\begin{equation*}
\mu_{n}(L)=\mathbb{P}(\{P_{n}\in L\})=\mathbb{P}(\{P_{n}|\Lambda_{m}^{(d)}\in\tilde{L}\})=\mathbb{E}\left(1_{\{P_{n}|\Lambda_{m}^{(d)}\in\tilde{L}\}}\right)
\end{equation*}
Hence
\begin{eqnarray*}
|\mu(L)-\mu_{n}(L)|&=&\\
\left|\mathbb{E}\left(1_{\{P|\Lambda_{m}^{(d)}\in\tilde{L}\}}-1_{\{P_{n}|\Lambda_{m}^{(d)}\in\tilde{L}\}}\right)\right|& \leq & \\
\mathbb{E}\left(\left|1_{\{P|\Lambda_{m}^{(d)}\in\tilde{L}\}}-1_{\{P_{n}|\Lambda_{m}^{(d)}\in\tilde{L}\}}\right|\right) & \leq & \\
\mathbb{P}\left( \mathcal{A}(\Lambda_{m}^{(d)})\not\subset\Lambda_{n}\right) & \leq & (2m+1)^{d}\cdot\frac{(2\nu+1)^{d\cdot\lceil\frac{n-m}{\nu}\rceil}}{\lceil\frac{n-m}{\nu}\rceil!}\\
 & = & o(e^{-\alpha n}),\ \forall\alpha>0
\end{eqnarray*}
 \hfill $\blacksquare$

\section{Rigorous Numerical Estimates and Final Remarks}
\label{conclusion}
We know from section\!~\ref{SLLNs} that
\[\lim_{i\rightarrow\infty}\rho_{i}=\rho=\mathbb{E}(P_{0})=\mathbb{P}\left(\{P_{0}=1\}\right)=1-\mathbb{P}\left(\{P_{0}=0\}\right)\]
Moreover, the construction carried out in section\!~\ref{construction} has shown that $P_{0}$ is determined by the relative order of a super-exponentially small\footnote{in the sense of inequality (\!~\ref{estimation}).} random number of the i.i.d. uniform random variables $\omega_{x},x\in\mathbb{Z}^{d}$. Therefore, bearing definitions~\ref{ppogs} and~\ref{thermodynamiclimit} in mind, it is possible to identify by inspection disjoint subsets of both $\{P_{0}=1\}$ and $\{P_{0}=0\}$ and hence determine lower and upper estimates for $\rho$. 

We clarify these ideas by means of a concrete example:

Set $d=1$, $\nu=1$ and $S_{\nu}^{(d)}=S_{1}^{(1)}=\left\{(0,0,0,)\right\}$ as in section~\ref{basic}. This set-up amounts to the one-dimensional version of the nearest neighbour exclusion scheme we touched on in section~\ref{introduction}.  Define $r$ to be the smallest non-negative integer where a local minimum with respect to $(\omega_{i})_{i\geq0}$ occur, that is $r(\omega)\stackrel{\mathrm{def}}{=}\inf\{i\in\mathbb{Z}_{+}:\omega_{i}<\omega_{i+1}\}$. Similarly set $l(\omega)\stackrel{\mathrm{def}}{=}-\sup\{i\in\mathbb{Z}_{-}:\omega_{i-1}>\omega_{i}\}$.

Now consider the partition $\{E_{ij};\ i,j\in\overline{\mathbb{Z}}_{+}\}$ of $\mathcal{B}_{(0,1)}^{\mathbb{Z}}$ defined by\footnote{$\overline{\mathbb{Z}}_{+}\stackrel{\mathrm{def}}{=}\mathbb{Z}_{+}\cup\{+\infty\}$.} 
\[E_{ij}\stackrel{\mathrm{def}}{=}\{\omega\in(0,1)^{\mathbb{Z}}:l(\omega)=i \mbox{ and }r(\omega)=j\}\]
and set $p_{ij}=\mathbb{P}(E_{ij})$.

It is easy to see that
\[\{P_{0}=1\}=\{\omega\in(0,1)^{\mathbb{Z}}:P(\omega)(0)=1\}=\biguplus_{i,j=0}^{\infty}E_{2i,2j}\]
Hence \[\rho=\mathbb{P}\left(\{P_{0}=1\}\right)=\sum_{i,j=0}^{\infty}\mathbb{P}(E_{2i,2j})=\sum_{i,j=0}^{\infty}p_{2i,2j}\]
Since \ $\mathbb{P}(E_{i,j})$\small
\begin{eqnarray*}
&=&\frac{i\cdot j}{(i+j+3)\cdot (i+1)!\cdot (j+1)!}+\frac{2}{(i+j+3)\cdot (i+j+2)\cdot (i+j+1)\cdot i!\cdot j!}\\
&+&\frac{1}{(i+j+3)\cdot (i+j+2)}\cdot \left(\frac{i}{(i+1)!\cdot j!}+\frac{j}{(j+1)!\cdot i!}\right)\ \ \ \ \ \ \ \ \ \ \ \ \ \forall i,j\in\mathbb{Z}_{+}
\end{eqnarray*} 
\normalsize
by elementary combinatorics. It is possible to generate lower bounds for $\rho$ by summing up the first \emph{even-even} probabilities and to generate upper bounds by summing up the first \emph{non-even-even} probabilities. Thus
\begin{eqnarray*}
&&0,4304=\sum_{i,j=0}^{2}p_{2i,2j}\leq\rho \\
&\leq & 1-\left(\sum_{0\leq i\leq 2 \atop 0\leq j \leq 1}p_{2i,2j+1}+ \sum_{0\leq i\leq 1 \atop 0\leq j \leq 2}p_{2i+1,2j}+\sum_{0\leq i,j \leq 1}p_{2i+1,2j+1}\right)=0,4339
\end{eqnarray*}
It is interesting to observe that these estimates are in agreement with the assymptotic first moment of the number of adsorpted particles calculated by J. K. Mackenzie~\cite{Mackenzie} by means of analytical tools such as Laplace transforms, differential equations, etc...\footnote{It is shown therein that $\mathbb{E}(|P_{n}|)\sim0,4324 \cdot 2n$ for the one-dimensional nearst neighbour exclusion scheme we touched on above, where $|P_{n}|$ is short for $\sum_{x\in\Lambda_{n}^{(d)}} P_{n}(x)$ as in section~\ref{introduction}. In fact, Mackenzie thought in terms of unit intervals instead of sites/particles and made calculations for the expected vacant space in the jamming limit instead of the expected number of occupied sites. Here we have translated his terminology to our current language.}

Of course the same readiness is not present in higher dimensional settings. Yet it is possible to obtain non trivial estimates for $\rho$ at least for simple models such as the two-dimensional nearest neighbour exclusion scheme and, to our knowledge, even in this case no rigorous estimates for $\rho$ exist.\vspace{5mm}

\noindent
\textbf{Parking Processes as Interacting Particle Systems}\footnote{We refer the interested reader to R. Durrett's review article~\cite{IPSs} for a quick, yet meaningful, introduction to this subject.}:\\
One might think of the parking processes $P_{n}$,\small$n\geq0$ \normalsize \footnote{As we already paid attention to, we are abusing terminology and using the jamming limit $P_{n}$ to denote the underlying parking process.} described in definition~\ref{parkingprocessonboxes} as discrete-time versions (i.e. skeletons) of ``finite'' interacting particle systems  $(\eta_{t}^{(n)})_{t\geq0}$
with initial condition $\eta_{0}=0^{\mathbb{Z}^{d}}$ and \emph{flip} rates $c(x,\eta)$ given by 
\[c(x,\eta)=1_{\{x\in\Lambda_{n}^{(d)}\}}\cdot1_{\{\theta_{x}(\eta)|\Lambda_{\nu}^{(d)}\in S_{\nu}^{(d)} \}} \ \ \ \ \ \forall x\in\mathbb{Z}^{d}, \eta\in\{0,1\}^{\mathbb{Z}^{d}}\]
It is not difficult to see that this is indeed the case, provided $S_{\nu}^{(d)}$ is a decreasing set \footnote{i.e. for $\xi_{1},\xi_{2}\in S_{\nu}^{(d)}$, if $\xi_{1}\preccurlyeq\xi_{2}$, then $\xi_{2}\in S_{\nu}^{(d)}\Rightarrow\xi_{1}\in S_{\nu}^{(d)}$; where `$\preccurlyeq$' is taken for the usual partial order of \emph{spin} configurations.}. (The foregoing correspondence being not correct otherwise, since  non-decreasing parking schemes generate non-markovian parking processes.)

In this setting it would be natural to define the properly infinite interacting particle system $(\eta_{t})_{t\geq0}$, defined by the \emph{flip} rates 
\[c(x,\eta)=1_{\{\theta_{x}(\eta)|\Lambda_{\nu}^{(d)}\in S_{\nu}^{(d)} \}} \ \ \ \ \ \forall x\in\mathbb{Z}^{d}, \eta\in\{0,1\}^{\mathbb{Z}^{d}}\]
and  same  initial condition $\eta_{0}=0^{\mathbb{Z}^{d}}$ and construct (i.e. couple) $(\eta_{t}^{(n)})_{t\geq0}$,\small $n\geq0$ and \normalsize $(\eta_{t})_{t\geq0}$ together making use of the same \emph{Poisson marks} to conclude that
\[\lim_{t\rightarrow\infty}\eta_{t}^{(n)}\stackrel{\mathrm{def}}{=} \eta_{\infty}^{(n)}\stackrel{(n)}{\longrightarrow}\eta_{\infty}\stackrel{\mathrm{def}}{=}\lim_{t\rightarrow\infty}\eta_{t}\ \ \ \mathrm{a.s.}\]
Since the sequence $(\eta_{\infty}^{(n)})_{n\geq0}$ is probabilisticaly indistinguishable from the the sequence of jamming limits $(P_{n})_{n\geq0}$ from the foregoing sections (i.e. $(\eta_{\infty}^{(n)})_{n\geq0}\sim(P_{n})_{n\geq0}$), we may interpret $\eta_{\infty}$ as $P$, ``the limit jamming limit'' of section~\ref{construction} (definition~\ref{thermodynamiclimit}).

All this reasoning is meaningful and in fact it proves in straghtfoward fashion the existence of the limit measure $\mu$ (at least for the markovian case). Yet one should bear in mind that the proofs carried out in sections~\ref{SLLNs} and~\ref{seb} rely on the specifics of section~\ref{construction}. Put it in other terms,  the sole existence of $\mu$ is not enough  to accomplish the results of this paper.\vspace{5mm}\\ 
\textbf{Parking Processes with Boundary Conditions Other than Null}:\\
It is easy to see that whatever boundary condition we use in definition~\ref{parkingprocessonboxes} we always obtain the same thermodynamic jamming limit $P$ of definition~\ref{thermodynamiclimit}. More precisely, if the sequence $(\tilde{P}_{n})_{n\geq1}$ is defined in the same probability space of the former sequence $({P}_{n})_{n\geq1}$; exactly as prescribed in definition~\ref{parkingprocessonboxes}, with the sole exception that we modify step(i) to
\[\mbox{\textbf{Step(\~{\i})}\hspace{15mm}}\ \ \tilde{P}_{n}|\Lambda_{n}^{(d)}\leftarrow0^{\Lambda_{n}^{(d)}};
\mbox{\hspace{10mm}} \tilde{P}_{n}|\left(\mathbb{Z}^{d}\setminus\Lambda_{n}^{(d)}\right)\leftarrow\eta|\left(\mathbb{Z}^{d}\setminus\Lambda_{n}^{(d)}\right)\]
for an arbitrary $\eta\in\{0,1\}^{\mathbb{Z}^{d}}$, then
\[\tilde{P}_{n}\longrightarrow P\longleftarrow P_{n} \ \ \ \mathrm{a.s.}\]
To see this we only have to modify lemma~\ref{inclusion} slightly as follows:\vspace{3mm}\\
\textit{For each $\omega \in \{0,1\}^{\mathbb{Z}^{d}}$ and $x\in\mathbb{Z}^{d}$, if $\Lambda_{n}\supset\mathcal{A}\left(\{x\}\right)(\omega)$, then $P_{n+\nu}$ and $P_{\mathcal{A} \left(\{x\}\right)}$ are identical over $\mathcal{A}\left(\{x\}\right)$; that is $\left.P_{n+\nu}\right|\mathcal{A}\left(\{x\}\right)=
\left.P_{\mathcal{A} \left(\{x\}\right)}\right|\mathcal{A}\left(\{x\}\right)$};\vspace{3mm}\\
and to take $n$ sufficiently large in theorem~\ref{convergence} so that $\mathcal{A}(\Lambda_{m}^{(d)})\subset \Lambda_{n-\nu}^{(d)}$.\vspace{10mm}\\
\large\textbf{Acknowledgements}:\normalsize\vspace{3mm}\\
The author is thankful to Geraldine G\'{o}es Bosco and Pablo Augusto Ferrari for useful conversations about this paper held at the Institute of Mathematics and Statistics of S\~{a}o Paulo State University (IME-USP) and to FAPESP (Funda\c{c}\~{a}o de Amparo \`{a} Pesquisa do Estado de S\~{a}o Paulo) for financial support.\vspace{0mm}


\end{document}